\magnification=1100
\hoffset= .22in
\voffset=-.3in
\input amstex
\documentstyle{amsppt}
\input xypic
\NoBlackBoxes
\NoRunningHeads
%%%%%%%%%%%%%%%%%%%%%%%%%%%%%%%%%%%%%%%%%%%%%%%%%%%%%%%

\define\wt{\widetilde}                                %
\define\ov{\overline} 

\define\prf{\demo{\underbar{Proof}}}
\define\endpf{\enddemo}
\define\dfn#1{\definition{\bf\underbar{Definition #1}}}

%macros for arrows

\define\surj{\twoheadrightarrow}                      %
                          %

%%%%%%%%%%%%%%%%%%%%%%%%%%%%%%%%%%%%%%%%%%%%%%%%
%\input bht.top.tex

\document
\baselineskip 20pt

\topmatter
\title
Subexponential group cohomology and the $K$-theory of Lafforgue's algebra $\Cal A_{max}(\pi)$
\endtitle
\vskip.2in
\author
R. Ji (IUPUI), C. Ogle (OSU)
\endauthor
\affil
Feb. 2004; revised April 2005
\endaffil
\email
ogle@math.mps.ohio-state.edu, ronji@math.iupui.edu
\endemail
\endtopmatter

%%%%%%%%%%%%%%%%%%%%%%%%%%%%%%%%%%%%%%%%
\newpage
\vskip.2in
\centerline{\bf\underbar{\S0\;\; Introduction}}
\vskip.3in

Let $K_*(X)$ denote the generalized homology of $X$ with coefficients in the
spectrum $\underline{\underline{K}}(\Bbb C)$  [A],  $\pi$ a finitely generated discrete group and $B\pi$ its classifying space. For any Banach algebra $A(\pi)$ with $\Bbb C[\pi]\subseteq A(\pi)$, there is an assembly map 
$K_*(B\pi) \to K^t_*(A(\pi))$ for the topological $K$-theory of $A(\pi)$. 
In this paper we are interested in the \lq\lq maximal unconditional completion\rq\rq of $\Bbb C[\pi]$ in the reduced group $C^*$-algebra $C^*_r(\pi)$. This algebra, denoted $\Cal A_{max}(\pi)$, was introduced by Lafforgue in [La] (the definition is recalled below).
The main result of the paper is

\proclaim{\bf\underbar{Theorem 1}} Let $\pi$ be a finitely-generated group and
$x= (x_n,x_{n-2},\dots)$ denote an element of $K_n(B\pi)\otimes\Bbb C\cong (H_*(B\pi)\otimes K_*(\Bbb C))\otimes\Bbb C$. Suppose there exists a cohomology class $[c]\in H^n(B\pi;\Bbb C)$ represented by an $n$-cocycle of subexponential growth with $<c,x_n>\ne 0$. Then $x$ is sent by the assembly map to a nonzero element of $K_n^t(\Cal A_{max}(\pi))\otimes\Bbb C$. In particular, if all rational homology classes of $\pi$ are detected by cocycles of this type, then the assembly map for the topological $K$-theory of $\Cal A_{max}(\pi)$ is rationally injective.
\endproclaim

An $n$-cochain $c$ is $\lambda$-exponential ($\lambda > 1$) if there is a constant $C$ (depending on $\lambda$)  for which $\parallel c(g_1,\dots,g_n)\parallel\le C(\lambda^{\sum_{i=1}^n L(g_i)})$ for all $[g_1,\dots,g_n]$ in $C_n(B\pi)$, where $L(_-)$ is the standard word-length function associated to a set of generators of $\pi$. It is subexponential if it is $\lambda$-exponential for all $\lambda > 1$. For finitely-generated groups two distinct presentations of $\pi$ yield linearly equivalent word-length functions, so the notion of $\lambda$-exponential does not depend on the choice of presentation of $\pi$. The proof of Theorem 1 is based on the recent work of Puschnigg [P], [P1], in which the theory of local cyclic homology is developed. Given the results of these two papers, our task essentially reduces to showing that subexponential group cocycles extend to continuous local cyclic homology cocycles on $\Cal A_{max}(\pi)$. This is accomplished using relatively standard techniques from cyclic homology and appealing to two key results from [P].
\vskip.2in

In an earlier draft of this paper, the authors proved the weaker result that cohomology classes of polynomial growth paired with $K^t_*(\Cal A_{max}(\pi))$ in the manner described by the above theorem, following the ideas of [CM]. It is clear that our original approach could have been extended to include a more restrictive subexponential growth condition modeled on the subexponential technical algebra described in [J1], [J2]. We are indebted to the referee for suggesting this stronger result as well as its connection to [P]. There is substantial evidence that the class of groups satisfying the property that all of its rational homology is detected by cocycles of polynomial growth is quite large, possibly containing all finitely-generated discrete groups whose first Dehn function is of polynomial type [O1, Conj. A]. If this pattern persists for Dehn functions of subexponential growth, one might expect a similar density of subexponential cohomology classes for Olshanskii groups.
\vskip.2in

The second author would like to thank D. Burghelea for many helpful conversations in the preparation of this paper.

%%%%%%%%%%%%%%%%%%%%%%%%%%%%%%%%%%%%%%%%

\newpage
\vskip.2in
\centerline{\bf\underbar{\S1\;\;  Proof of the main result}}
\vskip.3in

We recall some basic facts about Hochschild and cyclic
homology ([C1], [C2], [L]), consistent with the notation and
conventions of [P]. We assume throughout that the base field is
$\Bbb C$, with all algebras and tensor products being over $\Bbb
C$. For an associative unital algebra $A$, let $\Omega^n A = A\otimes
A^{\otimes n}$, and set $\Omega^*(A) = \underset n\to{\oplus}
\Omega^n(A)$. Define $b_n : \Omega^n(A)\to \Omega^{n-1}(A)$ by the
equation
$$
b_n(a_0,a_1,\dots, a_n) =
\sum_{i=0}^{n-1} (-1)^i (a_0,\dots, a_ia_{i+1},\dots, a_n)
+(-1)^n (a_na_0, a_1,\dots, a_{n-1})
$$
This is a differential, and the resulting complex
$(\Omega^*(A),b)$ is the Hochschild complex of $A$, also denoted
by $C_*(A)$. Its homology is denoted $HH_*(A)$. The Hochschild
complex admits a degree one chain map $B_*$ given by
$$
\split
B_n(a_0,a_1,\dots,a_n) = &\sum_{i=o}^n (-1)^{ni}(1, a_i, a_{i+1},\dots, a_n, a_0, a_1\dots, a_{i-1})\\
& -\sum_{i=0}^n(-1)^{ni}(a_i,1,a_{i+1},\dots,a_n,a_0,\dots,a_{i-1})
\endsplit
$$
and this defines a 2-periodic bicomplex $\Cal B(A)_{**}$ with
$\{\Cal B(A)_{2p,q} = \otimes^{q+1} A$, $\Cal B(A)_{2p-1,q} = 0$.
The differentials are $b_q:\Cal B(A)_{2p,q}\to \Cal B(A)_{2p,q-1}$
and $B_q:\Cal B(A)_{2p,q}\to \Cal B(A)_{2p-2,q+1}$. The total
complex $T_*(A) = (\underset{2p+q = *}\to\oplus \Cal B(A)_{2p,q}, b + B)$ of the resulting
bicomplex $\Cal B(A)_{**} = \{\Cal B(A)_{**}; b,B\}$ is clearly
$\Bbb Z/2$-graded (and contractible). There is
a decreasing filtration by subcomplexes $F^kT_*(A) =
(b(\Omega^k(A))\underset n\to\bigoplus \Omega^{n\ge k}(A), b +
B)$, and the periodic cyclic complex $CC_*^{per}(A)$ is the
completion of $T_*(A)$ with respect to this filtraton:
$CC_*^{per}(A) = \underset n\to\varprojlim T_*(A)/F^kT_*(A)$, with
homology $HC_*^{per}(A)$. For locally convex topological algebras,
projectively completing the tensor products in the complex
$\Omega^*(A)$ yields topological Hochschild and periodic cyclic
homology. \vskip.1in

\underbar{\bf\underbar{Note on notation}}: In [P] and [P1] the author
defines $\Omega^n(A)$ as $\Omega^n A = \wt A\otimes A^{\otimes n}$ where $\wt
A = A_+ = A\oplus\Bbb C$ denotes $A$ with a unit adjoined.  This is referred to as the Hochschild
complex of $A$, and its homology  as the Hochschild homology of
$A$. However, this definition actually produces the normalized
Hochschild complex $C_*(A_+)_n$ of $A_+$, and its homology is
$HH_*(A_+) = HH_*(A)\oplus HH_*(\Bbb C)$. Dividing out by by the
additional copy of $\Bbb C$ occuring in $\wt A$ in dimension $0$
of $C_*(A_+)_n$ yields the reduced complex $C_*(A_+)_{red}$, and
the inclusion $C_*(A)\hookrightarrow C_*(A_+)_{red}$ induces an
isomorphism in homology, but is not an isomorphism of complexes
[L, 1.4.2]. In order for the chain-level isomorphisms in [P, Lemma
3.2 and 3.3] to hold, one must use $C_*(\Bbb C[\pi])$, not the
larger $C_*(\Bbb C[\pi]_+)_{red}$. This decomposition leads to the
homogeneous decomposition stated in [P, Lemma 3.9], a result used
in this paper. For these reasons, we will assume throughout that our algebra $A$ is unital,
and that $\Omega^*(A)$ denotes the usual Hochschild complex of $A$.\vskip.2in

We recall the definition of  local cyclic homology [P], [P1]. Let
$A$ be a Banach algebra, $U$ its unit ball. Given a compact subset $K$ of $U$, $A_K$
is the completion of the subalgebra of $A$ generated by $K$  with respect to the largest
submultiplicative seminorm $\eta$ such that $\eta(K)\le 1$. Next, for an auxiliary Banach algebra $A'$,
let $\eta(_-)_{N,m}$ be the largest seminorm on $\Omega^*(A')$ satisfying
$$
\eta_{N,m}(a_0,a_1,a_2\dots ,a_n)_{N,m}\le\frac{1}{c(n)!}(2 + 2c(n))^m
N^{-c(n)} \parallel a_0\parallel_{A'}\cdot\dots\cdot \parallel
a_n\parallel_{A'} \tag1.1
$$
for each $n$, where $N\ge 1$, $m\in \Bbb N$ and $c(2n) = c(2n+1) = n$. For each $N$ the
boundary maps $b$ and $B$ are bounded with respect to the family of seminorms
$\{\eta(_-)_{N,m}\}_{m\in \Bbb N}$. Thus the completion of $T_*(A')$ with respect to
this family of seminorms again yields a complex, denoted $T_*(A')_{(N)}$. In particular
this applies when $A' = A_K$. Note that an inclusion of compact subspaces
$K\hookrightarrow K'$ induces a continuous homomorphism of Banach algebras
$A_K\to A_{K'}$, and that if $N < N'$ the identity map on $T^*(A')$ induces a natural
morphism of complexes
$T_*(A')_{(N)}\to T_*(A')_{(N')}$.

\proclaim{\bf\underbar{Definition 1.2 [P, Def. 3.4]}} In terms of
the above notation, the local cyclic homology of a Banach algebra
$A$ is
$$
HC_*^{loc}(A) = \underset{N\to \infty}\to{\underset {K\subset
U}\to\varinjlim}H_*(T_*(A_K)_{(N)})
$$
\endproclaim

A Banach space $X$ is said to have the \underbar{Grothendieck approximation property}
if the set of finite rank operators is dense in the space of
bounded operators $\Cal L(E)$ with respect to the compact-open
topology. A Banach algebra has this property when its underlying
Banach space does. In [P1] Puschnigg proves that for Banach
algebras satisfying this property, the above inductive system used
in the definition of $HC_*^{loc}(A)$ can be replaced by a much
smaller countable directed system. Precisely, let $V_0\subset
V_1\subset\dots\subset V_j\subset\dots\subset A$ be any increasing
sequence of finite-dimensional subspaces of $A$ for which $\bigcup
V_j$ is dense in $A$. Let $\ov B_r$ be the closed ball of radius
$r$ in $A$. Set $K_j = V_j\cap \ov B_{\frac{j}{j+1}}$.

\proclaim{\bf\underbar{Proposition 1.3 [P, Prop. 3.5; P1, Th.
3.2]}} If $A$ is a Banach algebra with the Grothendieck approximation
property, then
$$
HC^{loc}_*(A) = \underset{N\to \infty}\to{\underset
{j\to\infty}\to\varinjlim}H_*(T_*(A_{K_j})_{(N)})
$$
\endproclaim

We now specialize to the case $A$ is a Banach algebra completion
of $\Bbb C[\pi]$. There is a well-known decomposition of the
Hochschild complex ([B], [Ni]) $C_*(\Bbb C[\pi]) = \bigoplus
C_*(\Bbb C[\pi])_{<g>}$ the sum being indexed over the set
$\{<g>\}$ of conjugacy classes of $\pi$. The boundary maps $b_n$
preserve summands, so that for each $<g>$ $(C_*(\Bbb
C[\pi])_{<g>},b_*)$ is a subcomplex and direct summand of
$(C_*(\Bbb C[\pi]),b_*)$. Because the differential $B$ also
preserves this decomposition, it extends to the bicomplex $\Cal
B(A)_{**}$. These sum decompositions then yield sum-decompositions
of the corresponding algebraic homology groups
$$
HH_*(\Bbb C[\pi]) = \oplus HH_*(\Bbb C[\pi])_{<g>}
\tag1.4
$$
and for each $<g>$ a summand $HC_*^{per}(\Bbb C[\pi])_{<g>}$ of $HC_*^{per}(\Bbb C[\pi])$,
where $HH_*(\Bbb C[\pi])_{<g>}$ resp. $HC_*^{per}(\Bbb
C[\pi])_{<g>}$ is the homology of $(C_*(\Bbb C[\pi])_{<g>},b_*)$
resp. $(CC_*^{per}(\Bbb C[\pi])_{<g>},b_*)$. In general, extension
of these decompositions to the topological complexes associated to
completions of the group algebra is problematic. However, for
weighted $\ell^1$ completions, the summands persist (cf. [J1]).
For local cyclic homology, one has a similar result. To describe
it we need some more terminology.

Let $\ell^1(\pi)$ be the $\ell^1$ algebra of the discrete group
$\pi$. As in  [Bo], let $\ell^1_{\lambda}(\pi)$ denote the
completion of $\Bbb C[\pi]$ with respect to the largest seminorm
$\nu_{\lambda}$ satisfying $\nu_{\lambda}(g)\le \lambda^{L(g)}$
where $L$ is the the word length function on $\pi$ associated to a finite symmetric set of generators, and
$\lambda\ge 1$. Obviously $\ell^1_1(\pi) = \ell^1(\pi)$, and for
each $\lambda > 1$, $\ell^1_{\lambda}(\pi)$ is a subalgebra of
$\ell^1(\pi)$ and a Banach algebra, with
$\{\ell^1_{\lambda}(\pi)\}_{\lambda\ge 1}$ forming an inductive
system. The Banach algebra $\ell^1(\pi)$ has the Grothendieck
approximation property described above. Given a finite generating
set $S$ of $\pi$ with associated word-length metric $L_S$, let
$V_j$ be the linear span of $\{g\ |\ L_S(g)\le j\}$, with $K_j$
defined in terms of $V_j$ as before.

\proclaim{\bf\underbar{Proposition 1.5 [P, Lemma 3.7]}} For each
generating set $S$ and associated word-length function $L$ on
$\pi$ associated with $S$, there is a natural isomorphism
$$
"\underset{j\to\infty}\to\lim " \ell^1(\pi)_{K_j} \cong
"\underset{\lambda > 1}\to{\underset{\lambda\to 1}\to\lim}
"\ell^1_{\lambda}(\pi)
$$
\endproclaim

The quotation marks here indicate that these limits are actually
objects in the category of Ind-Banach algebras (i.e., the category
whose objects are formal inductive limits over the category of
Banach algebras -  [P1]). As a consequence, one has the following
analogue in local cyclic homology of the \lq\lq Principe
d'Oka\rq\rq due to Bost [Bo].

\proclaim{\bf\underbar{Theorem 1.6 [P, Cor. 3.8, Lemma 3.9]}} Let
$\pi$ be finitely-generated with generating set $S$ and associated
word-length function $L$. Then there is an isomorphism
$$
HC_*^{loc}(\ell^1(\pi))\cong
\underset{N\to\infty}\to{\underset{\lambda\to 1, \lambda >
1}\to\lim} H_*(CC_*(\ell^1_{\lambda}(\pi))_{(N)})
$$
Moreover $HC_*^{loc}(\ell^1(\pi))$ admits a decomposition as a
topological direct sum indexed by the set of conjugacy classes of
$\pi$.
\endproclaim

Next we recall the following definition and results, due to Lafforgue
[La]:  $\Cal A(\pi)$ is a sufficiently large good completion of
the group algebra $\Bbb C[\pi]$ inside of $C^*_r(\pi)$ if i) it is
admissible, ii) it is a Fr\'echet algebra with seminorms
$\{\eta_i(_-)\}$ satisfying $\eta_i(\sum_g \lambda_g g) =
\eta_i(\sum_g|\lambda_g| g)$; $|\lambda_g|\le |\lambda'_g|\forall
g$ implies $\eta_i(\sum_g\lambda_g g)\le \eta_i(\sum_g\lambda'_g g)$ and
iii) it contains the Banach subalgebras $\ell^1_{\lambda}(\pi)$
for all $\lambda > 1$. There is a maximal such subalgeba of
$C^*_r(\pi)$, denoted $\Cal A_{max}(\pi)$, which is the completion
of $\Bbb C[\pi]$ in the seminorm $\parallel f\parallel_{max} :=
\parallel |f| \parallel_r$, $\parallel_-\parallel_r$ being the
reduced $C^*$ norm. For our purposes, we will only state the
following result for $\Cal A_{max}(\pi)$; a more general result is
given in [P]. In what follows, $<e>$ denotes the conjugacy class of the identity element; the summand indexed by this element is referred to in [P] as the homogeneous summand (cf. also [B]).

\proclaim{\bf\underbar{Theorem 1.7 [P, Proposition 4.5], [CM]}} There is a decomposition
$$
HC_*^{loc}(\Cal A_{max}(\pi))\cong HC_*^{loc}(\Cal
A_{max}(\pi))_{<e>}\oplus (HC_*^{loc}(\Cal A_{max}(\pi))_{<e>})^C
$$
analogous to Burghelea's decomposition for $HC_*^{per}(\Bbb
C[\pi])$. Moreover the inclusions
$\ell^1_{\lambda}(\pi)\hookrightarrow \Cal A_{max}(\pi)$ induce an
isomorphism on the homogeneous summands of local cyclic homology groups
$$
HC_*^{loc}(\ell^1(\pi))_{<e>}\cong HC_*^{loc}(\Cal A_{max}(\pi))_{<e>}
$$
\endproclaim

We define a cohomology theory on $B\pi$ which we will pair with the right-hand side of (1.6).
Let $S$ and $L$ be as above. Let $B_*(\pi)$ denote the
non-homogeneous bar complex on $\pi$. Then $H^*(B\pi;\Bbb C)$ is
computed as the cohomology of the cocomplex $C^*(B\pi;\Bbb C) =
Hom(B_*(\pi);\Bbb C)$, a typical $n$-cochain represented by a map
$\phi:(\pi)^n\to\Bbb C$.

\proclaim{\bf\underbar{Definition 1.8}} A cochain
$\phi:(\pi)^n\to\Bbb C$ is called \underbar{$\lambda$-exponential}
if there exists a constant $C_{\lambda} > 0$ with
$$
|\phi(g_1,g_2,\dots,g_n)|\le C_{\lambda}\left(\lambda^{\sum_{i=1}^n L(g_i)}\right)
$$
for all $(g_1,g_2,\dots,g_n)\in B_n(\pi)$.
\endproclaim

For fixed $\lambda$, the set of $\lambda$-exponential cochains
form a subcocomplex $C_{\lambda}^*(B\pi;\Bbb C)$ of $C^*(B\pi;\Bbb
C)$, and the \underbar{$\lambda$-exponential cohomology} of $B\pi$
is $H_{\lambda}^*(B\pi;\Bbb C) := H^*(C_{\lambda}^*(B\pi;\Bbb
C))$. A cohomology class $[c]\in H^n(B\pi;\Bbb C)$ is
\underbar{$\lambda$-exponential} if it lies in the image of the
natural map $H_{\lambda}^*(B\pi;\Bbb C)\to H^*(B\pi;\Bbb C)$
induced by the inclusion $C^*_{\lambda}(B\pi;\Bbb
C)\hookrightarrow C^*(B\pi;\Bbb C)$. The \underbar{subexponential
cocomplex} is $C^*_{se}(B\pi;\Bbb C) := \underset\lambda >
1\to\bigcap C^*_{\lambda}(B\pi;\Bbb C)$, and the
\underbar{subexponential cohomology of $\pi$} is
$H^*_{se}(B\pi;\Bbb C) := H^*(C^*_{se}(B\pi;\Bbb C))$. Finally, a
cohomology class $[c]$ as above is subexponential if it is in the
image of the canonical map $H^*_{se}(B\pi;\Bbb C)\to H^*(B\pi;\Bbb
C)$. 

In order to describe the pairing, we construct a variant of Puschnigg's local theory. 
First, note that $HC^{loc}_*(A)$ is a $\Bbb Z/2$-graded theory whose algebraic
counterpart is not cyclic homology but periodic cyclic homology. Returning to the 
bicomplex $\Cal B(A)_{**}$, we let $\Cal {BC}(A)_{pq} = \{\Cal B(A)_{pq}\}$ for $p\ge 0$, 
and $0$ for $p < 0$. Algebraically this is the quotient of  $\Cal B(A)_{**}$ by the subbicomplex 
$\{\Cal B(A)_{pq}\}_{p < 0}$. Set $TC_*(A) = (\underset{2p+q = *}\to\oplus \Cal BC(A)_{2p,q}, b + B)$. This is a connective complex (usually denoted $CC_*(A)$) whose homology computes the algebraic cyclic homology of $A$. Let $C_*^{cy}(A) = \{C^{cy}_n(A) = C_n(A)/(1-\tau_{n+1}),b\}_{n\ge 0}$ denote the usual cyclic complex of $A$, formed as the quotient of the Hochschild complex by $(1-\tau_*)$, where $\tau_{n+1}(a_0,\dots,a_n) = (-1)^{n+1}(a_n,a_0,\dots,a_{n-1})$. There is a projection map $TC_*(A)$ which sends $\Cal{BC}_{0q}(A) = C_q(A)$ surjectively to $C^{cy}_q(A)$, $(1-\tau_{q+1})(C_q(A))\mapsto 0$, and sends $\Cal{BC}_{pq}(A)$ to $0$ for $p > 0$. As we are over a field of char. $0$,  this is a quasi-isomorphism. Now the completions used in the definition of local cyclic homology apply to both $TC_*(A)$ and $C^{cy}_*(A)$. Thus

\proclaim{\bf\underbar{Definition 1.9}} For a Banach algebra
$A$, let
$$
\gather
HCC_*^{iloc}(A) = \underset{N\to \infty}\to{\underset {K\subset
U}\to\varinjlim}H_*(TC_*(A_K)_{(N)})\\
HCC_*^{loc}(A) = \underset{N\to \infty}\to{\underset {K\subset
U}\to\varinjlim}H_*(C^{cy}_*(A_K)_{(N)})
\endgather
$$
\endproclaim
where the notation is as in Def. 1.2 above. The projection maps of complexes $T_*(A')\surj TC_*(A')\surj C^{cy}_*(A')$ are clearly bounded with respect to the seminorms used above for any of the auxiliary Banach algebras $A_K$, yielding homomorphisms on homology groups
$$
HC_*^{loc}(A)\to HCC_*^{iloc}(A)\to HCC_*^{loc}(A)
\tag1.10
$$
induced by the obvious projection maps of complexes.
Here $HCC_*^{loc}$ stands for \lq\lq \underbar{connective} local cyclic homology\rq\rq, as distinguished from $HC_*^{loc}$, while \lq\lq$iloc$\rq\rq refers to the \lq\lq intermediate local\rq\rq connective theory which simply serves as a bridge between the left and right hand sides of (1.10). Finally, we note that for algebras with the Grothendieck approximation property, one could similarly define homology groups using the smaller countable direct limits described above in Proposition 1.3. As we are interested in the case $A = \ell^1(\pi)$, we will use the directed systems of Theorem 1.6.

\proclaim{\bf\underbar{Definition 1.11}} Let
$\pi$ be finitely-generated with generating set $S$ and associated
word-length function $L$. Then
$$
HCC_*^{aploc}(\ell^1(\pi)) :=
\underset{N\to\infty}\to{\underset{\lambda\to 1, \lambda >
1}\to\lim} H_*(C^{cy}_*(\ell^1_{\lambda}(\pi))_{(N)})
$$
\endproclaim
By the same reasoning as in [P, Lemma 3.9], we see that
$HCC_*^{aploc}(\ell^1(\pi))$ admits a decomposition as a
topological direct sum indexed by the set of conjugacy classes of
$\pi$. Moreover, the canonical map $HC_*(\Bbb C[\pi])\to  HCC_*^{aploc}(\ell^1(\pi))$ obviously preserves the decomposition. In particular, $HC_*(\Bbb C[\pi])_{<e>}$ maps to $HCC_*^{aploc}(\ell^1(\pi))_{<e>}$. Note that, unlike the original local cyclic theory of [P], there is no reason to believe that $HCC_*^{aploc}(\ell^1(\pi))$ and $HCC_*^{loc}(\ell^1(\pi))$ (defined as $HCC_*^{loc}(_-)$ of an Ind-Banach algebra in the sense of [P1]) agree. Nevertheless, there is an evident map $HC_*^{loc}(\ell^1(\pi))\to HCC_*^{aploc}(\ell^1(\pi))$ which we may fit into a diagram
$$
\diagram
K_*^{top}(\Bbb C[\pi])\rrto\dto & & HC_*(\Bbb C[\pi])\dto\\
K_*^{top}(\ell^1(\pi))\rto\dto & HC_*^{loc}(\ell^1(\pi))\dto\rto & HCC_*^{aploc}(\ell^1(\pi))\dto|>>\tip\\
K_*^{top}(\Cal A_{max}(\pi))\rto & HC_*^{loc}(\Cal A_{max}(\pi))\rto & HCC_*^{aploc}(\ell^1(\pi))_{<e>}
\enddiagram
\tag1.12
$$
The vertical maps on the left are induced by the evident inclusion of topological algebras. The group $K_*^{top}(\Bbb C[\pi])$ denotes the topological $K$-theory of $\Bbb C[\pi]$ equipped with the fine topology, and the top horizontal map is the Chern character of [T]. The second and third horizontal maps are the Chern character of [P] and [P1] to local cyclic homology. The two maps to the group in the lower right corner are induced by Theorem 1.6, Theorem 1.7 and the evident decomposition of $HC_*^{aploc}(\ell^1(\pi))$ into summands indexed by conjugacy classes. The commutativity of the diagram follows from the compatibility of Chern characters. Theorem A now follows from

\proclaim{\bf\underbar{Lemma 1.13}} Let $[\phi]\in H^n_{se}(B\pi;\Bbb C)$. Let $[c]$ denote the image of $[\phi]$ in $H^n(B\pi;\Bbb C)$, and let $\tau_c$ denote the standard $n$-dim. cyclic group cocycle on $\Bbb C[\pi]$ formed by choosing a normalized representative $c$ of $[c]$ and extending it over $C^{cy}_n(\Bbb C[\pi])$. Then $\phi$ induces a map $(\phi)_*: HCC_n^{aploc}(\ell^1(\pi))_{<e>}\to\Bbb C$ such that the composition
$$
HC_n(\Bbb C[\pi])\to HCC_n^{aploc}(\ell^1(\pi))\surj HCC_n^{aploc}(\ell^1(\pi))_{<e>}\overset (\phi)_*\to\longrightarrow\Bbb C
$$
equals the map induced by $[\tau_c]$.
\endproclaim

\prf For fixed $N,m$ let $D_{N,m,n}= c(n)!(2 + 2c(n))^{-m}
N^{c(n)}$ (cf. (1.1)). Now also fix $\lambda > 1$, and consider an element $x = \sum_i\gamma_i (g_{0i}, g_{1i},\dots, g_{ni})$ in the completed complex $C^{cy}_n(\ell^1_{\lambda}(\pi))_{(N)}$. For $\tau_c$ as above, we have
$$
\gathered
| \tau_c(x)|\\
\le \sum\Sb i\\g_{0i}g_{1i}\cdot\dots\cdot g_{ni} = 1\endSb  |\gamma_i| | \tau_c(g_{0i},g_{1i},\dots,g_{ni}) |\\
= \sum_i |\gamma_i| | c([g_{1i},\dots ,g_{ni}])|\\
\le \sum_i |\gamma_i| C \lambda^{L(g_{1i})}\cdot\dots\cdot \lambda^{L(g_{ni})}\\
\le (C D_{N,m,n})\eta_{N,m}\left(\sum_i\gamma_i(g_{0i}, g_{1i}\dots, g_{ni})\right) = (C D_{N,m,n})\eta_{N,m}(x) < \infty
\endgathered
\tag1.14
$$
We are using the seminorm defined in (1.1) for the auxiliary algebra $A' = \ell^1_{\lambda}(\pi)$ with norm $\parallel \sum\gamma_i g_i\parallel_{\ell^1_{\lambda}(\pi)} := \sum_i |\gamma_i| \lambda^{L(g_i)}$.
Thus, for each $\lambda > 1$ and $N > 0$, $\tau_c$ extends to a continuous cyclic $n$-cocycle on 
$C^{cy}_*(\ell^1_{\lambda}(\pi))_{(N)}$, defining a homomorphism $H_n(C^{cy}_*(\ell^1_{\lambda}(\pi))_{(N)})\to\Bbb C$. These extensions are obviously compatible with the inclusions associated to the directed system occuring in the definition of $HCC_*^{aploc}(\ell^1(\pi))$ in (1.11) above, and factor through the projection $HCC_*^{aploc}(\ell^1(\pi))\surj HCC_*^{aploc}(\ell^1(\pi))_{<e>}$. Denoting the induced homomorphism on the direct limit by $(\phi)_*$ yields the result.\hfill // 

Combining this with the commuting diagram now completes the proof of main result stated in the introduction.

\endpf

%%%%%%%%%%%%%%%%%%%%%%%%%%%%%%%%%%%%%%%%

\enddocument